\documentclass[oneside,reqno]{amsart}

\usepackage{amssymb}
\usepackage{amsmath}
\usepackage{amsthm}
\usepackage{amsbsy}
\usepackage{bm}
\usepackage{hyperref}
\date{\today}
\usepackage{cite}
\usepackage{array}
\usepackage{xcolor}
\usepackage{tikz}
\usepackage{graphics}
\usepackage{subcaption}
\usepackage{bbm}

\makeatletter
\newcommand{\distas}[1]{\mathbin{\overset{#1}{\kern\z@\sim}}}

\newcommand{\distras}[1]{
	\mathbin{\overset{#1}{\kern\z@\resizebox{\wd\mybox}{\ht\mysim}{$\sim$}}}
}

\hypersetup{
	colorlinks=true,
	linkcolor=blue,
	filecolor=dviolet,
	urlcolor=blue,}

\theoremstyle{theorem}
    \newtheorem{theorem}{Theorem}
    \newtheorem{lemma}[theorem]{Lemma}
    
    \newtheorem{corollary}[theorem]{Corollary}

\theoremstyle{definition} 
    \newtheorem{definition}[theorem]{Definition}
    
    \newtheorem{result}[theorem]{Result}
    \newtheorem{remark}[theorem]{Remark}
    \newtheorem{example}[theorem]{Example}
    \newtheorem{exercise}[theorem]{Exercise}

\def\suchthat{\; : \;}

\def\Z{\mathbb{Z}}

\def\Z{\mathbb{Z}}

\def\v{{\bf v}}

\def\tends{\rightarrow}

\def\l{\left}
\def\r{\right}
\def\<{\langle}
\def\>{\rangle}

\newcommand{\E}{\mbox{E}}

\newcommand\Tr{{\mbox{Tr}}}

\newcommand\mnote[1]{} 
\newcommand\be{\begin{equation*}}

\newcommand\ee{\end{equation*}}

\newcommand\ben{\begin{equation}}
\newcommand\een{\end{equation}}
\newcommand\bes{\begin{eqnarray*}}
\newcommand\ees{\end{eqnarray*}}

\newcommand\bex{\begin{exercise}}
\newcommand\eex{\end{exercise}}
\newcommand\beg{\begin{example}}
\newcommand\eeg{\end{example}}
\newcommand\benu{\begin{enumerate}}
\newcommand\eenu{\end{enumerate}}
\newcommand\beit{\begin{itemize}}
\newcommand\eeit{\end{itemize}}
\newcommand\berk{\begin{remark}}
\newcommand\eerk{\end{remark}}
\newcommand\bdefn{\begin{defintion}}
\newcommand\edefn{\end{definition}}
\newcommand\bthm{\begin{theorem}}
\newcommand\ethm{\end{theorem}}
\newcommand\bprf{\begin{proof}}
\newcommand\eprf{\end{proof}}
\newcommand\blem{\begin{lemma}}
\newcommand\elem{\end{lemma}}

\newcommand{\var}{\mbox{\rm Var}}

\newcommand{\sm}{{\raise0.3ex\hbox{$\scriptstyle \setminus$}}}

\def\l{\left}
\def\r{\right}

\def\tends{\rightarrow}

\def\CHI{\mathchoice%
{\raise2pt\hbox{$\chi$}}%
{\raise2pt\hbox{$\chi$}}%
{\raise1.3pt\hbox{$\scriptstyle\chi$}}%
{\raise0.8pt\hbox{$\scriptscriptstyle\chi$}}}
\def\smalloplus{\raise1pt\hbox{$\,\scriptstyle \oplus\;$}}

\title[LSD of Toeplitz and Hankel matrices with dependent entries]{Limiting spectral distribution of Toeplitz and Hankel matrices with dependent entries}

\author{Shambhu Nath Maurya\\\\
	Statistics and Mathematics Unit,
	Indian Statistical Institute\\ Kolkata 700108, India}

\date{\today}
\thanks{shambhumath4 [at] gmail.com}

\begin{document}

\begin{abstract}
This article deals with the limiting spectral distributions (LSD) of symmetric Toeplitz and Hankel matrices with dependent entries. For any fixed positive integer $m$, we consider these $n \times n$ matrices with entries $\{Y^{(m)}_j / \sqrt n \ ;{j\geq 0}\}$, where $Y^{(m)}_j=  \sum_{r=-m}^{m} X_{j+r}$ and $\{X_k\}$ are i.i.d. with mean zero and variance one. We provide an explicit expression for the moment sequences of the LSDs. 
As a special case, this article provides an alternate proof for the LSDs of these matrices when the entries are i.i.d. with mean zero and variance one. The method is based on the moment method. The idea of proof can also be applied to other patterned random matrices, namely reverse circulant and symmetric circulant matrices.
\end{abstract}

\maketitle


\section{\large\textbf{Introduction and main results}}
Let $A_n$ be an $n\times n$ real symmetric random matrix. Then the \textit{empirical spectral distribution} (ESD) of $A_n$ is defined as
\begin{equation} \label{eqn:ESD}
F_{A_n} (x) = \frac{1}{n}\sum_{k=1}^{n} \mathbbm{1} (\lambda_k \leq x),
\end{equation}
where $\lambda_1,\lambda_2,\ldots, \lambda_n$ are eigenvalues of $A_n$. The weak limit of the ESD is known as the \textit{limiting spectral distribution} (LSD) of $A_n$, if it exists, in the probability or almost surely.

The LSD of random matrices is quite an interesting area of research in random matrix theory. 
 In 1958, Wigner \cite{wignerLSD1958} considered symmetric matrices with i.i.d. real Gaussian entries (Wigner matrix) and proved that the LSD is the \textit{semicircle law}. Bose and Mitra \cite{bose_mitra_LSD_Rc} established that the LSD of reverse circulant matrices is a \textit{symmetrized Rayleigh} distribution. For Toeplitz and Hankel matrices,
 Bryc et al. \cite{bryc_lsd_06} and Hammond, Miller \cite{hammond_miler_ToeplitzLSD} proved the existence of their LSDs. Later in 2008, Bose and Sen \cite{bose_sen_LSD_EJP} presented a unified approach to finding the limiting spectral distributions of symmetric patterned random matrices.
For results on the LSD of sample covariance matrices and other  matrices,  we refer the readers to \cite{pastur_LSDfirst_Smatrix}, \cite{bai_jack_LMRA_book} and \cite{bose_patterned}. 
 
  Our paper aims to study the LSDs of symmetric Toeplitz and Hankel matrices when their entries have some dependent structure and satisfy some moment assumptions. In literature, several results have been obtained for different random matrices when the entries have some types of dependent structure and moment assumptions. For the LSD results on the dependent entries, we refer the readers to \cite{chatterjee_LSD+lindeberg_06}, \cite{arijit+hazra_LSD+wigner_15} for Winger matrices; \cite{bose_LSD_depend_EJP} for circulant type matrices; \cite{costel+magda_LSD+SRM_16}, \cite{banna+magda_LSD+SRM_15}  for some symmetric random matrices and  \cite{gernot+philippe_ORMTbook_11}, \cite {oliver+eckhard_LSD+NRM_12}, \cite{bose_patterned} for other types of random matrices. There are several existing methods to study the LSD of random matrices. The method of moments (see \cite{bryc_lsd_06}, \cite{bose_sen_LSD_EJP}), Pastur’s fundamental technique of Stieltjes transforms (see \cite{pastur_Pcondtion_first}, \cite{chatterjee_LSD+lindeberg_06}) and Normal approximations (see \cite{bose_mitra_LSD_Rc}, \cite{bose_LSD_depend_EJP})  
  are some important methods to deal with the LSD of random matrices. 
 
 Now we define the matrices that we deal with.
 \vskip2pt
  \noindent \textbf{Toeplitz matrix:} An $n\times n$ Toeplitz matrix is defined as 
  \begin{align*}
  	T^{(a)}_n=\left(\begin{array}{ccccc}
  		x_{0} & x_{-1} & x_{-2} & \cdots & x_{1-n}\\
  		x_{1} & x_{0} & x_{-1} & \cdots & x_{2-n}\\
  		x_{2} & x_{1} & x_{0} & \cdots & x_{3-n}
  		\\ \vdots &\vdots &\vdots & \ddots & \vdots\\
  		x_{n-1}& x_{n-2}& x_{n-3}& \cdots & x_{0}
  	\end{array}  \right).
  \end{align*}
  Its $(i,\;j)$-th element is $x_{j-i}$. The \textit{symmetric Toeplitz matrix} has its $(i,\;j)$-th element as $x_{|j-i|}$ and we denote it by $T_n$.
 \vskip2pt 
\noindent \textbf{Hankel matrix:}  An $n\times n$  Hankel matrix is defined by its $(i,\;j)$-th element as $x_{j+i-1}$. If $D_n=(\delta_{i-1,n-j})_{n\times n}$ is the {\it backward identity} matrix  of dimension $n$, then $D_nT^{(a)}_n$ has a form of  Hankel matrix and its look like as
  \begin{align*}
  	H_n=\left(\begin{array}{cccccc}
  		x_{n-1}& x_{n-2}& x_{n-3}& \cdots & x_{0} \\
  		 x_{n-2}& x_{n-3}& x_{n-4}& \cdots & x_{-1} \\
  		x_{n-3} & x_{n-4} & x_{n-5} &\cdots  & x_{-2}\\
  		\vdots & \vdots & \vdots & \ddots &   \vdots\\
  		x_{0} & x_{-1} & x_{-2} & \cdots & x_{1-n}
  	\end{array}  \right).
  \end{align*}
  Observe that its $(i,\;j)$-th element is $x_{n-(j+i)+1}$. In this article, we always consider Hankel matrix of the form  $D_nT^{(a)}_n$.
  
  The following results are known for the LSD of $T_n$ and $H_n$.
  \begin{result} \label{res:LSD_rcscTH}
  	 \textbf{(Theorems 1.1, 1.2; \cite{bryc_lsd_06})} Suppose the entries of $T_{n}$ and $H_n$ are i.i.d. random variables with mean zero and variance one. Then almost surely,  the ESD of $T_{n}$ and $H_n$ converges weakly to some symmetric distribution that has unbounded support.
  	  \end{result}
  Note that in Result \ref{res:LSD_rcscTH}, the LSDs were studied for independent entries with some moment assumption. The authors used the method of moments to prove the results. Now in this article, we consider these matrices when the entries are the sum of finitely many independent copies, and we study their LSDs. Here, we generalize the LSD results of Result \ref{res:LSD_rcscTH} with respect to the entries of the matrices. In a particular case, our results will yield Result \ref{res:LSD_rcscTH}. 
  
Now we state our main results.
The following theorem provides the LSD of Toeplitz matrices.
\begin{theorem}\label{thm:toeESD}
	For integer $m \geq 0$, let $T_{n,m}$ be a symmetric Toeplitz matrix with entries $\{\frac{Y^{(m)}_j}{\sqrt n};{j\geq 0}\}$, where $Y^{(m)}_j= \sum_{r=-m}^{m} X_{j+r} $ and $\{X_k\}$ are i.i.d. with mean zero and variance one. Then the LSD of $T_{n,m}$ is a symmetric probability distribution $\mathcal{L}_{T,m}$ almost surely,
	where the moment sequence $\{\beta^{(T)}_h\}_{h \geq 1}$ of $\mathcal{L}_{T,m}$ is given by
	\begin{align}\label{eqn:moment_beta_T}
		\beta^{(T)}_h  & = 
		\left\{\begin{array}{ll} 	 
			\gamma_T(p) \times \# G^{(T)}_{2p,m}  & \text{if}\ h=2p, \\\
			0 & \text{if}\ h=2p+1	 	 
		\end{array}\right.		 
	\end{align}
	with $\gamma_T(p)$ and $\# G^{(T)}_{2p,m}$ are as given in (\ref{eqn:lim_2pmoment_T1}) and (\ref{eqn:card_D02pm_T}), respectively.
\end{theorem}
The following theorem provides the LSD of Hankel matrices.
\begin{theorem}\label{thm:hanESD}
	For integer $m \geq 0$, let $H_{n,m}$ be a Hankel matrix with entries $\{\frac{Y^{(m)}_j}{\sqrt n};{j\geq 0}\}$, where $Y^{(m)}_j= \sum_{r=-m}^{m} X_{j+r} $ and $\{X_k\}$ are i.i.d. with mean zero and variance one. Then the LSD of $H_{n,m}$ is a symmetric probability distribution $\mathcal{L}_{H,m}$ almost surely,
	where the moment sequence $\{\beta^{(H)}_h\}_{h \geq 1}$ of $\mathcal{L}_{H,m}$ is given by
	\begin{align} \label{eqn:moment_beta_H}
		\beta^{(H)}_h  & = 
		\left\{\begin{array}{ll} 	 
			\gamma_H(p) \times \# G^{(H)}_{2p,m}  & \text{if}\ h=2p, \\\
			0 & \text{if}\ h=2p+1	 	 
		\end{array}\right.		 
	\end{align}
	with $\gamma_H(p)$ and $\# G^{(H)}_{2p,m}$ are as given in (\ref{eqn:lim_2pmoment_H1}) and (\ref{eqn:card_D02pm_H}), respectively.
\end{theorem}
The following theorem says about a property of the LSDs of $T_{n,m}$ and $H_{n,m}$.
\begin{theorem} \label{thm:unbdd_LSD}
	The distributions $\mathcal{L}_{T,m}$ and $\mathcal{L}_{H,m}$ have unbounded support.
\end{theorem}
\begin{remark}
	\noindent (i) Observe from (\ref{eqn:card_D02pm_T}) and (\ref{eqn:card_D02pm_H}) that $\# G^{(T)}_{2p,0} =\# G^{(H)}_{2p,0}=1$. Therefore Theorem \ref{thm:toeESD} and Theorem \ref{thm:hanESD} give  Result \ref{res:LSD_rcscTH}.  Our combinatorial techniques are different than the combinatorics  used in \cite{bryc_lsd_06}.
	
	\noindent (ii)	The LSDs of reverse circulant and symmetric circulant matrices with input entries as stationary linear processes (two-sided moving average processes) were studied in \cite{bose_LSD_depend_EJP} in $L_2$ convergence. 
	If we consider the entries are finite two-sided moving average processes (as in Theorem \ref{thm:toeESD}), then using our idea, the results of \cite{bose_LSD_depend_EJP} can be derived for almost sure convergence.
\end{remark}
  
We use the method of moments  to establish our results. 
Note from (\ref{eqn:ESD}) that the $h$-th moment of the ESD of a real symmetric matrix $A_n$ can be written in terms of the trace of $(A_n)^h$.
Thus, if a closed form of the trace of $(A_n)^h$ is known, then one can exploit it to study the LSD of $A_n$.

Now we provide a brief outline of the rest of the manuscript. In Section \ref{pre_ESD_RcSc}, we state some basic notation and definitions which will arise throughout the article. We prove Theorem \ref{thm:toeESD} and Theorem \ref{thm:hanESD} in Section \ref{sec:ESD_Toe} and Section \ref{sec:ESD_Han}, respectively. Finally, in Section \ref{sec:propert_LSD+card_EG}, we prove Theorem \ref{thm:unbdd_LSD} and calculate the cardinality of some combinatorial objects.

 \section{\large \textbf{Preliminaries}}\label{pre_ESD_RcSc}
In this section, we state some basic notation and prove a result which is used to prove our main theorems. We first recall some standard notation and state a result.  Suppose $F_{A_n}$ is the ESD of a real symmetric matrix $A_n$. Then
\begin{equation}\label{eqn:moment_Tr_formula}
\tilde{E}_h(F_{A_n}) = \frac{1}{n} \sum_{k=1}^n (\lambda_k)^h = \frac{1}{n} \Tr(A_n)^h,
\end{equation}     
      where $\tilde{E}_h$ denotes the $h$-th moment with respect to the \textit{empirical spectral measure} of $A_n$ 
       and $\Tr$ denotes the trace of a matrix. The above expression is usually known as the \textit{trace-moment formula}.
      
The following result provides sufficient condition for the convergence of ESD 
of a given real symmetric random matrix $A_n$.
\begin{result} \textbf{(Lemma 1.2.4, \cite{bose_patterned})} \label{res:M_1234} 
 Suppose there exists a sequence $\{\beta_h\}_{h \geq 1}$ such that 
\begin{enumerate}
		\item [$(M_1)$] for every $h\in \mathbb{N}$, $\E[ \tilde{E}_h(F_{A_n})] \rightarrow \beta_h$ as $n \to \infty$,
		\item [$(M_2)$] there is a unique distribution $F$ whose moment sequence is $\{\beta_h\}_{h \geq 1}$, and
		\item [$(M_3)$] for every $h\in \mathbb{N}$, $\sum_{n=1}^{\infty}\E \big[ \tilde{E}_h(F_{A_n})-\E[ \tilde{E}_h(F_{A_n})] \big]^4<\infty$.
	\end{enumerate}
	Then almost surely, $\{F_{A_n}\}$ converges weakly to $F$.
\end{result}

Now we start with some notions and results which will be used to prove our theorems. Let $\Omega$ be the space of all probability distributions with the finite second moment. It is well-known that the space $\Omega$ is a metric space with $W_2$-metric, which is defined as
$$W_2(F,G) = \Big[ \inf_{(X \distas{D} F, Y \distas{D} G)} \E(X-Y)^2 \Big]^{\frac{1}{2}}  \ \ \mbox{ for  } F,G \in \Omega,$$
where $X \distas{D} F$ denotes that  $X$ and $F$ have the same distribution. 
The following result provides a relation between weak convergence and convergence in  $W_2$-metric.
\begin{result} \textbf{(Lemma 1.3.1, \cite{bose_patterned})} \label{res:W2+weak_relation}
	The metric $W_2$ is complete and for $F_n,F \in \Omega$, $W_2(F_n, F) \to 0$ if and only if $F_n$ converges weakly to $F$ and $\E(F_n)^2 \to \E(F)^2$.
\end{result}
The following result provides a link between the trace of matrices and $W_2$ distance of the ESDs of matrices.
\begin{result}\textbf{ (Lemma 1.3.2, \cite{bose_patterned})} \label{res:W2+trace_relation}
	Suppose $A, B$ are $n \times n$ real symmetric matrices. Then
	\begin{equation*} 
		[W_2(F_A,F_B)]^2 \leq \frac{1}{n} \Tr(A-B)^2,
	\end{equation*}
	where $F_A$ and $F_B$ are ESDs of $A$ and $B$, respectively.
\end{result}
The following corollary is a nice application of Result \ref{res:W2+weak_relation} and Result \ref{res:W2+trace_relation}.
\begin{corollary}\label{lem:LSD_A=A_0}
	Let $A_n$ be any symmetric random matrix with $(i,j)$-th entry as $(a_{i,j} / \sqrt n)$, where $\{a_{i,j}\}_{i,j\geq 0}$ are i.i.d. with mean zero and variance one. Furthermore, suppose $A_{n,0}$ is a matrix defined as 
	\begin{equation*} \label{eqn:A_0}
		(A_{n,0})_{i,j}	
		=  \left\{\begin{array}{ll} 
			(a_{i,j} / \sqrt n) & \text{if }   \ \     i \neq j, \\
			0 & \text{if } \ \    i=j,	 
		\end{array}\right. 
	\end{equation*}
	where $(A_{n,0})_{i,j}$ denotes the $(i,j)$-th entry of $A_{n,0}$.
	Then the LSD of $A_n$ is same as the LSD of $A_{n,0}$.
\end{corollary}
\begin{proof} 
	Note from Result \ref{res:W2+trace_relation} that
	\begin{equation*} 
		[W_2(F_{A_n},F_{A_{n,0}})]^2 \leq \frac{1}{n} \Tr(A_n-A_{n,0})^2=  \frac{1}{n} \sum_{i=1}^{n} (\frac{a_{i,i}}{\sqrt{n}})^2 \tends 0 \mbox{ as } n \tends \infty.
	\end{equation*}
	Combining the above convergence with Result \ref{res:W2+weak_relation} shows that the LSD of $A_n$ is same as the LSD of $A_{n,0}$.
\end{proof}

Now we  provide a theorem which gives a relation between the LSDs of matrices when the entries of matrices are independent and dependent ( the sum of finitely many independent copy), respectively.
\begin{theorem}   \label{thm:iid_uniform_ESD}
	Let $A_{n,m}$ be symmetric Toeplitz or Hankel matrix with entries $\{\frac{Y^{(m)}_j}{\sqrt n};{j\geq 0}\}$. Suppose the ESD of $A_{n,m}$ converges weakly almost surely to some fixed non-random distribution $F$ when $\{X_j\}$ satisfy the following assumption:
\begin{equation}\label{eqn:condition_ESD}
	\{X_j\} \mbox{ are independent}, \ \E(X_j)=0, \ \E(X_j^2)=1 \ \forall \ j \in \Z \mbox{ and }  \sup_{j \in \Z}\E(|X_j|^k)=\alpha_k<\infty \ \mbox{for}\ k\geq 3.
\end{equation}
	Then the same limit continues to hold if $\{X_j\}$ are i.i.d. with mean zero and variance one.
\end{theorem}
  
 \begin{proof} 
 	We prove this theorem for symmetric Toeplitz matrix only. For Hankel matrix, the arguments will be similar.
 	Let $T_{n,m}$ be a symmetric Toeplitz matrix with entries $\{Y^{(m)}_j/\sqrt n \ ;{j\geq 0}\}$, where $Y^{(m)}_j= \sum_{r=-m}^{m} X_{j+r}$ and $\{X_i\}$ are i.i.d. with mean zero and variance one.
 	Now for a fixed positive constant $K$, we define a set of new variables 
 	\begin{equation*} 
 		X'_i= \frac{1}{\sigma_{x,K}} \big[ X_i \mathbb{I}(|X_i|\leq K)- \E\big\{X_i \mathbb{I}(|X_i|\leq K)\big\} \big],
 	\end{equation*}
 	where $\sigma^2_{x,K} = \var[X_i \mathbb{I}(|X_i|\leq K)]$. Note that $\{X'_i\}$ are independent with mean zero and variance one. Since $\sigma^2_{x,K} \to 1$ as $K\to \infty$, $\{X'_i\}$ is bounded.
 	Let $T^{(K)}_{n,m}$ be Toeplitz matrix with entries $\{Y^{'(m)}_j /\sqrt n \ ;{j\geq 0}\}$, where $Y^{'(m)}_j= \sum_{r=-m}^{m} X'_{j+r}$.
 	Observe that $\{X'_i\}$ satisfy condition (\ref{eqn:condition_ESD}). Therefore, from the assumption of Theorem \ref{thm:iid_uniform_ESD}, the LSD of $T^{(K)}_{n,m}$ exists.
 	
 	Now note from Result \ref{res:W2+trace_relation} that
 	\begin{align} \label{eqn:S12esd}
 		\big[W_2 \big(F_{T_{n,m}}, F_{T^{(K)}_{n,m}} \big)\big]^2 
 		\leq \frac{1}{n} \Tr \big[T_{n,m} - T^{(K)}_{n,m} \big]^2 
 		&\leq  \frac{1}{n} \times 2n \sum_{j=0}^{n-1} \big(\frac{Y^{(m)}_j}{\sqrt n} - \frac{Y^{'(m)}_j}{\sqrt n} \big)^2 \nonumber \\
 		& =\frac{2}{n} \sum_{j=0}^{n-1} \big[\sum_{r=-m}^{m} (X_{j+r}- X'_{j+r}) \big]^2 \nonumber \\
 		& \leq 2 \sum_{r=-m}^{m} \big[ \frac{1}{n} \sum_{j=0}^{n-1} (X_{j+r}- X'_{j+r})^2 \big] \nonumber \\
 		& \tends 2\sum_{r=-m}^{m} \big[ \E (X_{j+r}- X'_{j+r})^2 \big],
 	\end{align}
 as $n \to \infty$. Now if we let $K \tends \infty$, then it is easy to show by an application of \textit{Dominated Convergence Theorem} that for each fixed $r$, $ \E (X_{j+r}- X'_{j+r})^2 \tends 0$ and hence $\big[W_2 \big(F_{T_{n,m}}, F_{T^{(K)}_{n,m}} \big)\big]^2 \tends 0.$
 Finally, using Result \ref{res:W2+weak_relation}, we get that the LSD of $T_{n,m}$
 is same as the LSD of $T^{(K)}_{n,m}$.
 This completes the proof
  of Theorem \ref{thm:iid_uniform_ESD}.
 \end{proof}

\section{\large \textbf{LSD of $T_{n,m}$}} \label{sec:ESD_Toe}
In this section, we prove Theorem \ref{thm:toeESD}.
First, we introduce and review some basic combinatorial
concepts, which will appear during the moment's calculation of the
LSD of $T_{n,m}$ and again in Section \ref{sec:ESD_Han} to calculate the LSD of $H_{n,m}$.

\begin{definition} \label{def:pair+epsilion} Let the set $[n]=\{1,2,\ldots,n\}$.
	
\noindent (i)  We call $\pi=\{V_{1},\ldots,V_{r}\}$ a partition of $[n]$ if the
		blocks  $V_{j} \,(1\leq j \leq r)$ are pairwise disjoint, non-empty
		subsets of $[n]$ such that $[n]=V_{1}\cup\cdots\cup V_{r}$. 
		
		\noindent (ii) The set of all partitions of $[n]$ is denoted by
		$\mathcal{P}(n)$, and the subset consisting of all pair partitions,
		that is, $\# V_{j}=2$, for all $1\leq j \leq r$, is denoted by
		$\mathcal{P}_{2}(n)$.

		\noindent (iii) The subset of $\mathcal{P}_{2 }(n)$ consisting of such pair partitions that each contains exactly one even number and one odd number is denoted by $\mathcal{P}^{oe}_{2}(n)$.
		
	\noindent (iv) Without loss of generality,
		any $\pi \in \mathcal P_2(2k)$ will be written as 
		$$\pi=(r_1,s_1)\cdots (r_k,s_k),\ \ \text{where} \ \ r_t<s_t,  \ \text{for} \ \ t=1,\ldots, k.$$ 
		Then, for $t=1,\ldots, k$, we can define
		the projection $\pi$ as
		\begin{align*} 
			\mbox{ $\pi(r_t)=\pi(s_t)= r_t$ and }	\epsilon_\pi(\ell)=\l\{\begin{array}{rll}
				1 & \mbox{ if } & \ell=r_t,\\ 
				-1& \mbox{ if } & \ell=s_t.
			\end{array} \r.
		\end{align*}
\end{definition}
We also define the following: 
$$	\chi_{[1,n]}(z) = 1\;\;\mbox{if $z \in [1,n]$  and zero otherwise}.$$

Now we recall a  convenient trace formula for Toeplitz matrices from \cite{liu_wang2011} which will be required to prove Theorem \ref{thm:toeESD}.
\begin{lemma} 
	\label{lem:traceTn_copy}
	Suppose $T_n$ is the Toeplitz matrices with entries $\{x_j;{j\geq 1}\}$. Then 
	\begin{equation}\label{trace formula T_n}
		\Tr[ (T_n)^h]
		= \sum_{i=1}^{n}\,\sum_{J_h \in A_h} x_{j_1} x_{j_2} \cdots x_{j_h}\, I(i, J_h),
	\end{equation}  
	where $J_{h}= (j_1, j_2, \ldots, j_{h}), I(i, J_h) = \prod_{k=1}^h \chi_{[1,n]}(i+\sum_{\ell=1}^{k} j_\ell )$ and
	\begin{align}\label{def:A_2p_TESD}
		A_{h} &=\big\{(j_1, j_2, \ldots, j_{h})\in \mathbb Z^{h}\suchthat \sum_{k=1}^{h} j_k=0, -n \le j_1,\ldots, j_{h}\le n\big\}.
	\end{align}
\end{lemma}
For the proof of the above lemma, we refer the reader to Lemma 2.4 of \cite{liu_wang2011}. 
Now we define some notions on vectors.
For a vector $J= (j_1,j_2, \ldots, j_p) \in \Z^p$, we define the multi-set $S_{|J|}$  as 
\begin{equation*} 
	S_{|J|}=\{ |j_1|, |j_2|, \ldots , |j_p| \}.
\end{equation*}

\begin{definition}\label{def:opposite_sign_Z}
	Suppose $v=(j_1, j_2, \ldots, j_p)$ is a vector in $\mathbb Z^p$. 
	
\noindent (i)  We say an entry $j_k$ of $v$ is \textit{pair matched} if $j_k$ appears exactly two times in the multi-set $\{j_1, j_2, \ldots, j_p\}$. Similarly, we can define a {\it triple matched} element if it appears three times. We say $v$ is pair matched if all its entries are pair matched.
	
	\noindent (ii) Two elements $j_k, j_\ell$ of $v$ are said to be {\it $\pm$ pair matched} if $j_k$ and $j_\ell$  are of opposite sign and $|j_k|$ appears exactly twice in $S_{|v|}.$
	For example, in both $(3,5,8,-5)$ and $(3,-5,8,5)$, the two entries $5$ and $-5$ are {\it $\pm$ pair matched} whereas in $(3,5,8,5)$, the two entries $5$ are not {\it $\pm$ pair matched} but pair matched. 
	For even $p$, we say that a vector $v$ is {\it $\pm$ pair matched}, if all its entries are {\it $\pm$ pair matched}. 
\end{definition}
The following remark is an observation which will be used later to prove theorems.

\begin{remark} \label{rem:pairmatch_obser}
	For a constant $c$, suppose $J_{p}=(j_1,j_2,\ldots, j_{p}) \in \Z^{p}$, $\sum_{k=1}^{p} j_k= c$ and the entries $\{j_1,j_2,\ldots, j_{p} \}$ are at least $\pm$ pair matched. Then for $p$ even, the maximum number of free entries in $J_{p}$ will be $p/2$ if and only if $J_{p}$ is {\it $\pm$ pair matched} and $c=0$. For $p$ odd, the maximum number of free entries in $J_{p}$ could be $\frac{p-1}{2}-1$, where $(-1)$ arises due to the constraint on the entries of $J_p$.
\end{remark}

Now with the help of the above result and trace formula, we prove Theorem \ref{thm:toeESD}.
\begin{proof}[Proof of Theorem \ref{thm:toeESD}] 
	In view of Theorem \ref{thm:iid_uniform_ESD}, it is enough to assume (\ref{eqn:condition_ESD}) for the entries $\{X_j\}$ of $T_{n,m}$ and prove the required almost sure convergence.
	Recall from Result \ref{res:M_1234} that to prove Theorem \ref{thm:toeESD}, it is sufficient to verify $(M_1), (M_2)$ and $(M_3)$ for $T_{n,m}$. 
	First we verify $(M_1)$.
	Note from (\ref{eqn:moment_Tr_formula}) and (\ref{trace formula T_n}) that
	\begin{align} \label{eqn:E_X_i_T}
		\E\big[ \tilde{E}_{h}(F_{T_{n,m}}) \big] 
		&=  \frac{1}{n^{\frac{h}{2}+1}} \sum_{i=1}^{n} \sum_{J_{h} \in A_{h}} \E[ Y^{(m)}_{j_1} Y^{(m)}_{j_2} \cdots Y^{(m)}_{j_{h}}] I(i, J_h) \nonumber \\
		& =  \frac{1}{n^{\frac{h}{2}+1}} \sum_{i=1}^{n} \sum_{J_{h} \in A_{h}} \E \Big[ \big(\sum_{r=-m}^{m}X_{j_1+r}\big) \big(\sum_{r=-m}^{m}X_{j_2+r}\big) \cdots \big(\sum_{r=-m}^{m}X_{j_{h}+r}\big) \Big] I(i, J_h) \nonumber \\
		&  =  \frac{1}{n^{\frac{h}{2}+1}}  \sum_{D_{h} \in G_{h,m}} \sum_{i=1}^{n} \sum_{J_{h} \in A_{h}}  \E[ X_{j_1+d_1} X_{j_2+d_2} \cdots X_{j_{h}+d_{h}}] I(i, J_h),
	\end{align} 
	where $J_{2p}= (j_1, j_2, \ldots, j_{h}), \ A_{h}$ is as in (\ref{def:A_2p_TESD}) and $G_{h,m}$ is defined as
	\begin{align} \label{eqn:E_2p,m_T}
		G_{h,m}= \big\{(d_1, d_2, \ldots, d_{h}) \in \mathbb{Z}^{h} : d_s \in \{-m, -(m-1), \ldots, -1,0,1, \ldots,m\} \ \forall \ s=1,2, \ldots, h \big\}.
	\end{align}
	Now we calculate the order of convergence of a typical term of (\ref{eqn:E_X_i_T})  for a fixed $(d_1, d_2, \ldots, d_{h}) \in G_{h,m}$. Note  that 
	\begin{equation} \label{eqn:T(p,m)_T}
		\sum_{J_{h} \in A_{h}}  \E[ X_{j_1+d_1} X_{j_2+d_2} \cdots X_{j_{h}+d_{h}}] = \sum_{K_{h} \in A^*_{h}}  \E[ X_{k_1} X_{k_2} \cdots X_{k_{h}}] = T(h,d), \mbox{ say},
	\end{equation}
	where for $K_{h}= (k_1, k_2, \ldots, k_{h})$,
	\begin{align}\label{def:A*_2p_TESD}
		A^*_{h} &=\big\{(k_1, k_2, \ldots, k_{h})\in \mathbb \Z^{h}\suchthat \sum_{s=1}^{h} k_s=  \sum_{s=1}^{h} d_s, -n \le k_s \le n \ \forall \ s=1,2, \ldots, h \big\}.
	\end{align}
	Since $\E[X_j] =0$ and $X_j=X_{-j}$ for each $j$, $ \E[ X_{k_1} X_{k_2} \cdots X_{k_{h}}]$ will be non-zero only when the entries $\{k_1, k_2, \ldots, k_{h}\}$ are at least $\pm$ pair matched. Thus, $ T(h,d) \leq O(n^\frac{h}{2})$. Observe from Remark \ref{rem:pairmatch_obser} that if $h$ is odd, then $ T(h,d)$ can have maximum order of convergence $O(n^{\frac{h-1}{2}})$. For $h$ even, say $2p$, $ T(p,d) = O(n^p)$ is only possible when $(k_1, k_2, \ldots, k_{2p})$ is exactly $\pm$ pair matched. Also note from Remark \ref{rem:pairmatch_obser} that a $\pm$ pair matched vector $(k_1, k_2, \ldots, k_{2p}) \in A^*_{2p}$ has $p$ degree of freedom if and only if   $\sum_{r=1}^{2p} d_r=0$. Hence 
	\begin{align}  \label{eqn:E_bound_rEven_T}
		T(h,d)
		&= \left\{\begin{array}{ll} 
			O(n^{\frac{h}{2}}) & \text{if} \   h \mbox{ is even, }  (k_1, k_2, \ldots, k_{h}) \mbox{ is $\pm$ pair matched and } \sum_{r=1}^{h}  d_r=0, \\
			o(n^{\frac{h}{2}}) & \text{otherwise}. 	 
		\end{array}\right. 
	\end{align}
	
	Note from (\ref{def:A_2p_TESD}) and (\ref{def:A*_2p_TESD}) that $A^*_{2p}= A_{2p}$ if $\sum_{r=1}^{2p}  d_r=0 $. Now if we define
	\begin{align} \label{eqn:D_2p,m,0_T}
		G^{(T)}_{2p,m} := \big\{ (d_1, d_2, \ldots, d_{2p}) \in G_{2p,m} : \sum_{r=1}^{2p} d_r=0 \big\},
	\end{align}
	then using (\ref{eqn:T(p,m)_T}) and (\ref{eqn:E_bound_rEven_T}) in (\ref{eqn:E_X_i_T}), we get
	\begin{align} \label{eqn:lim_ESD_evpT}
		&\lim_{n\to\infty} \E\big[ \tilde{E}_{h}(F_{T_{n,m}}) \big]  \nonumber \\ 
			& = \left\{\begin{array}{ll} 
			\displaystyle  \lim_{n\to\infty}  \frac{1}{n^p}  \sum_{D_{2p} \in G^{(T)}_{2p,m}} \sum_{i=1}^{n} \sum_{ \substack{{ K_{2p} \in A_{2p}} \\ {K_{2p} \mbox{ is $\pm$ pair matched}}}}  \E[ X_{k_1} X_{k_2} \cdots X_{k_{2p}}] I(i, K_{2p}) & \text{if }       h=2p, \\
			0 & \text{if } h=2p+1,
		\end{array}\right. \nonumber \\
		&=  \left\{\begin{array}{ll} 
			\displaystyle   \# G^{(T)}_{2p,m} \times \lim_{n\to\infty}  \frac{1}{n^p} \sum_{i=1}^{n} \sum_{ \substack{{ K_{2p} \in A_{2p}} \\ {K_{2p} \mbox{ is $\pm$ pair matched}}}}  \E[ X_{k_1} X_{k_2} \cdots X_{k_{2p}}] I(i, K_{2p}) & \text{if }       h=2p, \\
			0 & \text{if } h=2p+1,
		\end{array}\right.
	\end{align}
	where $\# G^{(T)}_{2p,m}$ is as in (\ref{eqn:card_D02pm_T}).

	Now we calculate the limiting term of (\ref{eqn:lim_ESD_evpT}).  Note from the trace formula (\ref{trace formula T_n}) that
	\begin{align*}
&\lim_{n\to\infty}  \frac{1}{n^p} \sum_{i=1}^{n} \sum_{ \substack{{ K_{2p} \in A_{2p}} \\ {K_{2p} \mbox{ is $\pm$ pair matched}}}}  \E[ X_{k_1} X_{k_2} \cdots X_{k_{2p}}] I(i, K_{2p}) \nonumber \\
&= \lim_{n\to\infty}  \frac{1}{n^p} \sum_{i=1}^{n}  \sum_{\pi \in \mathcal{P}_2(2p)} \sum_{k_1, \ldots, k_p=-n}^{n} \prod_{s=1}^{2p} \chi_{[1,n]} \big(i+\sum_{\ell=1}^{s} \epsilon_\pi(\ell) k_{\pi(\ell)} \big),
	\end{align*}
	where $\mathcal{P}_2(2p)$ is the set of all pair partition of $2p$ and $\epsilon_\pi, \pi(\ell)$ are as defined in Definition \ref{def:pair+epsilion}. Note that for each fixed $\pi \in \mathcal{P}_2(2p)$, 
	\begin{align*}
		\frac{1}{n^p} \sum_{i=1}^{n} \sum_{k_1, \ldots, k_p=-n}^{n} \prod_{s=1}^{2p} \chi_{[1,n]} \big(i+\sum_{\ell=1}^{s} \epsilon_\pi(\ell) k_{\pi(\ell)} \big) = 	\frac{1}{n^p} \sum_{i=1}^{n} \sum_{k_1, \ldots, k_p=-n}^{n} \prod_{s=1}^{2p} \chi_{[\frac{1}{n}, 1]} \big(\frac{i}{n}+\sum_{\ell=1}^{s} \epsilon_\pi(\ell) \frac{k_{\pi(\ell)}}{n}  \big),
	\end{align*}
	which can be considered as the Riemann sum of the following integral
	\begin{align*}
\int_{[0,1] \times [-1,1]^p} \prod_{s=1}^{2p} \chi_{[0,1]} \big(z_0+\sum_{\ell=1}^{s} \epsilon_\pi(\ell) z_{\pi(\ell)} \big)	\prod_{s=1}^{p} dz_s.
	\end{align*}
Hence we have
\begin{align}\label{eqn:lim_2pmoment_T1}
&\lim_{n\to\infty}  \frac{1}{n^p} \sum_{i=1}^{n} \sum_{ \substack{{ K_{2p} \in A_{2p}} \\ {K_{2p} \mbox{ is $\pm$ pair matched}}}}  \E[ X_{k_1} X_{k_2} \cdots X_{k_{2p}}] I(i, K_{2p}) \nonumber \\
&=  \sum_{\pi \in \mathcal{P}_2(2p)} \int_{[0,1] \times [-1,1]^p} \prod_{s=1}^{2p} \chi_{[0,1]} \big(z_0+\sum_{\ell=1}^{s} \epsilon_\pi(\ell) z_{\pi(\ell)} \big)	\prod_{s=1}^{p} dz_s \nonumber \\
& = \gamma_T(p), \mbox{ say}.	
\end{align}	
Using (\ref{eqn:lim_2pmoment_T1}) in (\ref{eqn:lim_ESD_evpT}), we get
	\begin{align} \label{eqn:lim_ESD_Tm}
	\lim_{n\to\infty} \E\big[ \tilde{E}_{h}(F_{T_{n,m}}) \big]   
			& =  \left\{\begin{array}{ll} 
			 \# G^{(T)}_{2p,m} \times \gamma_T(p)
			  & \text{if }       h=2p, \\
				0 & \text{if } h=2p+1,	 
			\end{array}\right. \nonumber \\
			&= \beta^{(T)}_h, \mbox{ say}.	 
\end{align}
	
	Now we verify condition $(M_2)$. Since $\{\beta^{(T)}_h\}$ is a limit of a moment sequence, we have that the sequence $\{\beta^{(T)}_h\}$ of (\ref{eqn:lim_ESD_Tm}) is a moment sequence of a measure, say $\mathcal{L}_{T,m}$ on $\mathbb{R}$.  For the uniqueness of $\mathcal{L}_{T,m}$, we use Riesz’s condition (Lemma 1.2.2 of \cite{bose_patterned}). Observe from (\ref{eqn:E_2p,m_T}) and (\ref{eqn:D_2p,m,0_T}) that $\# G^{(T)}_{2p,m} \leq (2m+1)^{2p}$ also note from (\ref{eqn:lim_2pmoment_T1}) that
	\begin{equation} \label{eqn:gammaT_leq}
		\gamma_T(p) \leq 2^p \times \# \mathcal{P}_2(2p)= 2^p \times \frac{(2p)!}{2^p p!}.
	\end{equation}
	 Thus
	\begin{align*}
		\liminf_{p \tends \infty} \frac{1}{p}	(\beta^{(T)}_{2p})^{\frac{1}{2p}}  \leq \liminf_{p \tends \infty} \frac{1}{p}\times 2 \big[\frac{(2p)!}{2^p p!} \big]^{\frac{1}{2p}} (2m+1) < \infty.
	\end{align*}
	Hence $\mathcal{L}_{T,m}$ is a unique distribution corresponding to the sequence $\{\beta^{(T)}_h\}$.
	
  The condition $(M_3)$ for $T_{n,m}$ can be verified by showing that
	\begin{equation*} 
		\E \big[ \tilde{E}_h(F_{T_{n,m}})-\E[ \tilde{E}_h(F_{T_{n,m}})] \big]^4 = O(n^{-2}) \ \mbox{ for all } h \in \mathbb{N},
	\end{equation*}
 For the proof of the above expression, see the proof of equation (2.6) of \cite{liu_wang2011}. The same arguments carry over here with some technical changes. We omit the details. 
 This completes the proof Theorem \ref{thm:toeESD}.
\end{proof}

\section{\large \textbf{LSD of $H_{n,m}$}} \label{sec:ESD_Han}
In this section, we prove Theorem \ref{thm:hanESD}.
First we start with the following trace formula for Hankel matrices.
\begin{lemma} 
	\label{lem:traceHn_copy}
	Suppose $H_n$ is the Hankel matrices with entries $\{x_j;{j\geq 1}\}$. Then 
	\begin{equation}\label{trace formula H_n}
		\Tr(H_n)^{p}=\left\{ \begin{array}{ll}
			\displaystyle  \sum_{i=1}^{n} \sum_{J_h \in C_h} x_{j_1} x_{j_2} \cdots x_{j_h}\,  I(i, J_h)
				 &  \hspace{-9pt}\mbox{ if $h$ is even,}
			\\	
			\displaystyle \sum_{i=1}^{n} \sum_{J_h \in C_{h,i}} x_{j_1} x_{j_2} \cdots x_{j_h} I(i, J_h)\, & \hspace{-7pt}\mbox{if $h$ is odd,}
		\end{array}
		\right.
	\end{equation}  
	where $J_{h}= (j_1, j_2, \ldots, j_{h}), I(i, J_h) = \prod_{\ell=1}^{p}\chi_{[1,n]}(i-\sum_{t=1}^{\ell}(-1)^t j_t)$ and for $i=1,2, \ldots,n$,
	\begin{align}\label{def:A_2p_HESD}
		C_{h} &=\big\{(j_1, j_2, \ldots, j_{h})\in \mathbb Z^{h}\suchthat \sum_{k=1}^{h} (-1)^k j_k=0, \ -n \le j_1,\ldots, j_{h}\le n\big\}, \nonumber \\
			C_{h,i} &=\big\{(j_1, j_2, \ldots, j_{h})\in \mathbb Z^{h}\suchthat \sum_{k=1}^{h} (-1)^k j_k=2i-1-n, \ -n \le j_1,\ldots, j_{h}\le n\big\}.\textbf{}
	\end{align}
\end{lemma}
We refer the reader to Lemma 2.5 of \cite{liu_wang2011} for the proof of the above lemma. 
Now we define a notion of {\it $\pm$ odd-even pair matched} vectors for the set of vectors from  $\mathbb Z^p$. First recall the notion of pair matched and triple matched elements from Definition \ref{def:opposite_sign_Z}.
\begin{definition}\label{def:odd-even_H}
	Suppose $v=(j_1, j_2, \ldots, j_p)$ is a vector in $\mathbb Z^p$. Two elements $j_k, j_\ell$ of $v$ are said to be {\it $\pm$ odd-even pair matched},  if $j_k=j_\ell$, $j_k$ appears exactly twice in $v$, once each at an odd and an even position. 
	For example; in $(5,5,8,3)$, the two entries $5$; in $(8,-5,-5,3)$, the two entries $-5$ are {\it $\pm$ odd-even pair matched}.
	Also, in the vector $(5,5,8,3,-5,-5)$,  the two entries $5$ and $-5$ are {\it $\pm$ odd-even pair matched}
	 whereas in $(3,5,8,5)$ and $(5,8,5,3)$, the two entries $5$ are not {\it $\pm$ odd-even pair matched}. 
	 
	 Note that, here we consider $j_k$ and $-j_k$ as two different entries. For example; in $(3,5,-5,8)$, the entry $5$ or $-5$ is not {\it $\pm$ odd-even pair matched}.
	 
	For even $p$, we say that a vector $v$ is {\it $\pm$ odd-even pair matched}, if all its entries are {\it $\pm$ odd-even pair matched}. 	For example; $(5,-5,8,8,-5,5)$, $(3,3,-6,-6)$ and $(-3,-6,-6,-3)$  are {\it $\pm$ odd-even pair matched} vectors. But  $(5,-5,8,8,5,-5)$, $(3,6,-6,-6)$ and $(-3,6,-6,-3)$  are not {\it $\pm$ odd-even pair matched} vectors.
\end{definition}
The following remark is an observation of the Definition \ref{def:odd-even_H}.
\begin{remark} \label{rem:pairmatch+-_obser}
	For fixed constant $c$, suppose $J_{p}=(j_1,j_2,\ldots, j_{p}) \in \Z^{p}$, $\sum_{k=1}^{p} (-1)^kj_k= c$ and the entries $\{j_1,j_2,\ldots, j_{p} \}$ are at least pair matched. Then for $p$ even, the maximum number of free entries in $J_{p}$ will be $p/2$ if and only if $J_{p}$ is {\it $\pm$ odd-even pair matched} and $c=0$. For $p$ odd, the maximum number of free entries in $J_{p}$ could be $\frac{p-1}{2}-1$, where $(-1)$ arises due to the constraint on the entries of $J_p$.
\end{remark}
Now using the above result and trace formula, we prove Theorem \ref{thm:hanESD}.
\begin{proof}[Proof of Theorem \ref{thm:hanESD}] 
	In view of Theorem \ref{thm:iid_uniform_ESD}, it is enough to assume (\ref{eqn:condition_ESD}) for the entries $\{X_j\}$ and prove the required almost sure convergence.
	First recall from Result \ref{res:M_1234} that to prove Theorem \ref{thm:hanESD}, it is sufficient to verify $(M_1), (M_2)$ and $(M_3)$ for $H_{n,m}$. 
	
	First we verify $(M_1)$.
	Note from (\ref{trace formula H_n}) that for odd and even values of $h$ we have different trace formulae for $(H_{n})^h$. So we calculate the limit of $\E[ \tilde{E}_h(F_{H_{n,m}})]$ for odd and even values of $h$ separately.
	First suppose $h$ is even, say $h=2p$. Then from (\ref{eqn:moment_Tr_formula}) and (\ref{trace formula H_n}), we get
	\begin{align} \label{eqn:E_X_i_H}
		\E\big[ \tilde{E}_{2p}(F_{H_{n,m}}) \big] 
		&=  \frac{1}{n^p} \sum_{i=1}^{n} \sum_{J_{2p} \in C_{2p}} \E[ Y^{(m)}_{j_1} Y^{(m)}_{j_2} \cdots Y^{(m)}_{j_{2p}}] I(i, J_{2p}) \nonumber \\
		& =  \frac{1}{n^p} \sum_{i=1}^{n} \sum_{J_{2p} \in C_{2p}} \E \Big[ \big(\sum_{r=-m}^{m}X_{j_1+r}\big) \big(\sum_{r=-m}^{m}X_{j_2+r}\big) \cdots \big(\sum_{r=-m}^{m}X_{j_{2p}+r}\big) \Big] I(i, J_{2p}) \nonumber \\
		&  =  \frac{1}{n^p}  \sum_{D_{2p} \in G_{2p,m}} \sum_{i=1}^{n} \sum_{J_{2p} \in C_{2p}}  \E[ X_{j_1+d_1} X_{j_2+d_2} \cdots X_{j_{2p}+d_{2p}}] I(i, J_{2p}),
	\end{align} 
	where $J_{2p}= (j_1, j_2, \ldots, j_{2p}), \ C_{2p}$ is as in (\ref{def:A_2p_HESD}) and $G_{h,m}$ is defined as in (\ref{eqn:E_2p,m_T}).
	Now we calculate the order of convergence of a typical term of (\ref{eqn:E_X_i_H})  for a fixed $(d_1, d_2, \ldots, d_{2p}) \in G_{2p,m}$. Note  that 
	\begin{equation} \label{eqn:T(p,m)_H}
		\sum_{J_{2p} \in C_{2p}}  \E[ X_{j_1+d_1} X_{j_2+d_2} \cdots X_{j_{2p}+d_{2p}}] = \sum_{K_{2p} \in C^*_{2p}}  \E[ X_{k_1} X_{k_2} \cdots X_{k_{2p}}] = T'(p,d), \mbox{ say},
	\end{equation}
	where for $K_{2p}= (k_1, k_2, \ldots, k_{2p})$,
	\begin{align}\label{def:A*_2p_HESD}
		C^*_{2p} &=\big\{(k_1, k_2, \ldots, k_{2p})\in \mathbb Z^{2p}\suchthat \sum_{s=1}^{2p}(-1)^s k_s=  \sum_{s=1}^{2p}(-1)^s d_s \mbox{ (mod $n$)}, 1\le  k_{s}\le  n   \ \forall \ s=1,2, \ldots, 2p \big\}.
	\end{align}
	Note from assumption (\ref{eqn:condition_ESD}) that $\E[X_j] =0$ for each $j$, therefore $ \E[ X_{k_1} X_{k_2} \cdots X_{k_{2p}}]$ will be non-zero only when the entries $\{k_1, k_2, \ldots, k_{2p}\}$ are at least pair matched. Thus, $ T'(p,d) \leq O(n^p)$. In fact, $T'(p,d) = O(n^p)$ is only possible when $(k_1, k_2, \ldots, k_{2p})$ is exactly pair matched. Observe from Remark \ref{rem:pairmatch+-_obser} that a pair matched vector $(k_1, k_2, \ldots, k_{2p}) \in C^*_{2p}$ has $p$ degree of freedom if and only if  $(k_1, k_2, \ldots, k_{2p})$ is $\pm$ odd-even pair matched and $\sum_{r=1}^{2p} (-1)^r d_r=0$. Hence 
	\begin{align}  \label{eqn:E_bound_rEven_H}
		T'(p,d)
		&= \left\{\begin{array}{ll} 
			O(n^p) & \text{if} \    (k_1, k_2, \ldots, k_{2p}) \mbox{ is $\pm$ odd-even pair matched and } \sum_{r=1}^{2p} (-1)^r d_r=0, \\
			o(n^p) & \text{otherwise}. 	 
		\end{array}\right. 
	\end{align}
	
	Note from (\ref{def:A_2p_HESD}) and (\ref{def:A*_2p_HESD}) that $C^*_{2p}= C_{2p}$ if $\sum_{r=1}^{2p} (-1)^r d_r=0 $. Now if we define
	\begin{align} \label{eqn:D_2p,m,0_H}
		G^{(H)}_{2p,m} := \big\{ (d_1, d_2, \ldots, d_{2p}) \in G_{2p,m} : \sum_{r=1}^{2p} (-1)^r d_r=0 \big\},
	\end{align}
	then using (\ref{eqn:T(p,m)_H}) and (\ref{eqn:E_bound_rEven_H}) in (\ref{eqn:E_X_i_H}), we get
	\begin{align} \label{eqn:lim_ESD_evpH}
		\lim_{n\to\infty} \E\big[ \tilde{E}_{2p}(F_{H_{n,m}}) \big] &= \lim_{n\to\infty}  \frac{1}{n^p}  \sum_{D_{2p} \in G^{(H)}_{2p,m}} \sum_{i=1}^{n} \sum_{ \substack{{ K_{2p} \in C_{2p}} \\ {K_{2p} \mbox{  is $\pm$ odd-even pair matched}}}}   \hskip-30pt \E[ X_{k_1} X_{k_2} \cdots X_{k_{2p}}] I(i,K_{2p}) \nonumber \\
		&= \# G^{(H)}_{2p,m} \times \lim_{n\to\infty}  \frac{1}{n^p} \sum_{i=1}^{n} \sum_{ \substack{{ K_{2p} \in C_{2p}} \\ {K_{2p} \mbox{  is $\pm$ odd-even pair matched}}}} \hskip-30pt \E[ X_{k_1} X_{k_2} \cdots X_{k_{2p}}] I(i,K_{2p}),
	\end{align}
	where $\# G^{(H)}_{2p,m}$ is as in (\ref{eqn:card_D02pm_H}).

	Now we calculate the limiting term of (\ref{eqn:lim_ESD_evpH}).  Note from the trace formula (\ref{trace formula H_n}) that
	\begin{align*}
		&\lim_{n\to\infty}  \frac{1}{n^p} \sum_{i=1}^{n} \sum_{ \substack{{ K_{2p} \in C_{2p}} \\ {K_{2p} \mbox{ is $\pm$ odd-even pair matched}}}}  \E[ X_{k_1} X_{k_2} \cdots X_{k_{2p}}] I(i, K_{2p}) \nonumber \\
		&= \lim_{n\to\infty}  \frac{1}{n^p} \sum_{i=1}^{n}  \sum_{\pi \in \mathcal{P}^{oe}_2(2p)} \sum_{k_1, \ldots, k_p=-n}^{n} \prod_{s=1}^{2p} \chi_{[1,n]} \big(i-\sum_{\ell=1}^{s} (-1)^\ell \epsilon_\pi(\ell) k_{\pi(\ell)} \big), \\
		&=	 \sum_{\pi \in \mathcal{P}^{oe}_2(2p)} \lim_{n\to\infty}   \frac{1}{n^p} \sum_{i=1}^{n}\sum_{k_1, \ldots, k_p=-n}^{n} \prod_{s=1}^{2p} \chi_{[\frac{1}{n}, 1]} \big(\frac{i}{n} - \sum_{\ell=1}^{s} (-1)^\ell \epsilon_\pi(\ell) \frac{k_{\pi(\ell)}}{n}  \big),
	\end{align*}
	where $\mathcal{P}^{oe}_2(2p)$, $\epsilon_\pi$ and $\pi(\ell)$ are as in Definition \ref{def:pair+epsilion}.
	Note that for each fixed $\pi \in \mathcal{P}^{oe}_2(2p)$, 
 the typical term in the above expression converge to  the following Riemann  integral
	\begin{align*}
		\int_{[0,1] \times [-1,1]^p} \prod_{s=1}^{2p} \chi_{[0,1]} \big(z_0- \sum_{\ell=1}^{s} (-1)^\ell \epsilon_\pi(\ell) z_{\pi(\ell)} \big)	\prod_{s=1}^{p} dz_s.
	\end{align*}
	Hence we have
	\begin{align}\label{eqn:lim_2pmoment_H1}
		&\lim_{n\to\infty}  \frac{1}{n^p} \sum_{i=1}^{n} \sum_{ \substack{{ K_{2p} \in C_{2p}} \\ {K_{2p} \mbox{ is $\pm$ odd-even pair matched}}}}  \E[ X_{k_1} X_{k_2} \cdots X_{k_{2p}}] I(i, K_{2p}) \nonumber \\
		&=  \sum_{\pi \in \mathcal{P}^{oe}_2(2p)} \int_{[0,1] \times [-1,1]^p} \prod_{s=1}^{2p} \chi_{[0,1]} \big(z_0-\sum_{\ell=1}^{s} (-1)^\ell \epsilon_\pi(\ell) z_{\pi(\ell)} \big)	\prod_{s=1}^{p} dz_s \nonumber \\
		& = \gamma_H(p), \mbox{ say}.	
	\end{align}	
	Using (\ref{eqn:lim_2pmoment_H1}) in (\ref{eqn:lim_ESD_evpH}), we get
	\begin{align} \label{eqn:lim_ESD_hev_Hm}
		\lim_{n\to\infty} \E\big[ \tilde{E}_{2p}(F_{H_{n,m}}) \big]   
		& = 
			\# G^{(H)}_{2p,m} \times \gamma_H(p).	 
	\end{align}

	Now suppose $h$ is odd, say $h=2p+1$. Then from (\ref{eqn:moment_Tr_formula}) and the trace formula (\ref{trace formula H_n}), we get
	\begin{align} \label{eqn:E_X_i_H_odd}
		\E\big[ \tilde{E}_{2p+1}(F_{H_{n,m}}) \big] 
		& = \frac{1}{n^\frac{2p+3}{2}} \sum_{i=1}^n \sum_{J_{2p+1} \in C_{2p+1, i}} \E[ Y^{(m)}_{j_1} Y^{(m)}_{j_2} \cdots Y^{(m)}_{j_{2p+1}}] \nonumber \\
		& =  \frac{1}{n^\frac{2p+3}{2}} \sum_{i=1}^n \sum_{J_{2p+1} \in C_{2p+1, i}} \E \Big[ \big(\sum_{r=0}^{m}X_{j_1+r}\big) \big(\sum_{r=0}^{m}X_{j_2+r}\big) \cdots \big(\sum_{r=0}^{m}X_{j_{2p+1}+r}\big) \Big] \nonumber \\
		&  =  \frac{1}{n^\frac{2p+3}{2}} \hskip-3pt \sum_{D_{2p+1} \in G_{2p+1,m}} \sum_{i=1}^n \sum_{J_{2p+1} \in C_{2p+1, i}}  \E[ X_{j_1+d_1} X_{j_2+d_2} \cdots X_{j_{2p+1}+d_{2p+1}}],
	\end{align} 
	where  $J_{2p+1}= (j_1, j_2, \ldots, j_{2p+1}), D_{2p+1}= (d_1, \ldots, d_{2p+1})$, $C_{2p+1,i}$ as in (\ref{def:A_2p_HESD}) and $G_{2p+1,m}$ as in (\ref{eqn:E_2p,m_T}).
	
	Since $\E[X_j] =0$, for a fixed $(d_1, d_2, \ldots, d_{2p+1}) \in G_{2p+1,m}$, $\E[ X_{j_1+d_1} X_{j_2+d_2} \cdots X_{j_{2p+1}+d_{2p+1}}]$ will be non-zero only when the entries $\{j_1+d_1, j_2+d_2, \ldots, j_{2p+1}+d_{2p+1}\}$ are at least pair matched. 
	Thus a typical term of (\ref{eqn:E_X_i_H_odd}) has the maximum contribution if one entry is triple matched and the remaining entries are pair-matched. In this case, the contribution will be of the order $O(n^{\frac{2p}{2} -1})$, 
	where $(-1)$ arises due to the constraint of $C_{2p+1,i}$, $\sum_{k=1}^{2p+1} (-1)^k j_k=2i-1 \mbox{ (mod $n$)}$. 
	So, we get
	\begin{align*} 
		\E\big[ \tilde{E}_{2p+1}(F_{H_{n,m}}) \big] 
		= \frac{1}{n^\frac{2p+3}{2}} O(n)  O(n^{p-1})= o(1).
	\end{align*}
	Hence from the above expression and (\ref{eqn:lim_ESD_hev_Hm}), we get
	\begin{align*}
		\beta^{(H)}_h := \lim_{n\to\infty}  \E\big[ \tilde{E}_{h}(F_{H_{n,m}}) \big] 
		= 
		\left\{\begin{array}{ll} 
		\# G^{(H)}_{2p,m} \times \gamma_H(p) & \text{if }       h=2p, \\
			0 & \text{if } h=2p+1. 	 
		\end{array}\right. 
	\end{align*}
	
	Now we verify condition $(M_2)$. First note that the sequence $\{\beta^{(H)}_h\}$ is a moment sequence of a measure, say $\mathcal{L}_{H,m}$ on $\mathbb{R}$, as $\{\beta^{(H)}_h\}$ is a limit of a moment sequence.  For the uniqueness of $\mathcal{L}_{H,m}$, we use Riesz’s condition. Note from  (\ref{eqn:D_2p,m,0_H}) that $\# G^{(H)}_{2p,m} \leq (2m+1)^{2p}$ and also $\gamma_H(p)\leq 2^p \times \frac{(2p)!}{2^p p!}$ (see (\ref{eqn:gammaT_leq})). Thus
	\begin{align*}
		\liminf_{p \tends \infty} \frac{1}{p}	(\beta^{(H)}_{2p})^{\frac{1}{2p}}  \leq \liminf_{p \tends \infty} \frac{1}{p}  \times \Big(\frac{(2p)!}{2^p p!} \Big)^{\frac{1}{2p}} (2m+1) < \infty.
	\end{align*}
	Hence $\mathcal{L}_{H,m}$ is a unique distribution corresponding to the sequence $\{\beta^{(H)}_h\}$.
	
 The condition $(M_3)$ for $H_{n,m}$ can be verified by showing that
\begin{equation*} 
	\E \big[ \tilde{E}_h(F_{H_{n,m}})-\E[ \tilde{E}_h(F_{H_{n,m}})] \big]^4 = O(n^{-2}) \ \mbox{ for all } h \in \mathbb{N},
\end{equation*}
which again follows from the arguments similar to those used to establish (2.6) of \cite{liu_wang2011}. Here we skip its proof. This completes the proof of Theorem \ref{thm:hanESD}.
\end{proof}

\section{\large \textbf{Properties of LSDs and Cardinality of  $G^{(H)}_{2p,m}$, $G^{(H)}_{2p,m}$}} \label{sec:propert_LSD+card_EG}
First we find the cardinalities of the sets. 

\noindent \textbf{Cardinality of $G^{(T)}_{2p,m}$:} First recall $G^{(T)}_{2p,m}$ from (\ref{eqn:D_2p,m,0_T}) that 
	\begin{align*} 
	G^{(T)}_{2p,m} := \big\{ (d_1, d_2, \ldots, d_{2p}) \in \mathbb{Z}^{2p} : d_r \in \{-m, -(m-1), \ldots, -1,0,1, \ldots,m\} \mbox{ and } \sum_{r=1}^{2p} d_r=0 \big\}.
\end{align*}
We use generating function technique to find the cardinality of $G^{(T)}_{2p,m}$. First observe that for fixed positive integers $p$ and $m$, $\# G^{(T)}_{2p,m}$ is the coefficient of $x^0$ in the expression of $(x^{-m}+ x^{-(m-1)}+ \cdots+ x^{-1}+x^0+x^1+\cdots+x^{m})^{2p}= f(x)$, say. Note that
\begin{align*}
	f(x) =   x^{-2mp}  (1+x+\cdots+ x^{2m})^{2p} &= x^{-2mp} (1-x^{2m+1})^{2p}  (1-x)^{-2p} \nonumber\\
	&  = x^{-2mp} \Big\{\sum_{k=0}^{2p}(-1)^k\binom{2p}{k}x^{(2m+1)k} \Big\}  \Big\{\sum_{s=0}^{\infty} \binom{2p+s-1}{2p-1}x^s \Big\}  \nonumber\\
	&= x^{-2mp} g(x), \mbox{ say}.
\end{align*}
Observe from the above expression that $\# G^{(T)}_{2p,m}$ is the coefficient of $x^{2mp}$ in $g(x)$ and therefore
\begin{align}  \label{eqn:card_D02pm_T}
	\# G^{(T)}_{2p,m}= 
		\sum_{\ell=0}^{\lfloor \frac{2mp}{2m+1} \rfloor} (-1)^\ell \binom{2p}{\ell} \binom{2p+2pm-(2m+1)\ell-1}{2p-1}. 
\end{align}
\vskip2pt

\noindent \textbf{Cardinality of $G^{(H)}_{2p,m}$:} First recall $G^{(H)}_{2p,m}$ from (\ref{eqn:D_2p,m,0_H}) that 	
\begin{align*}
	G^{(H)}_{2p,m} := \big\{ (d_1, d_2, \ldots, d_{2p})\in \mathbb{Z}^{2p} : d_r \in \{-m, -(m-1), \ldots, -1,0,1, \ldots,m\} \mbox{ and } \sum_{r=1}^{2p} (-1)^rd_r=0 \big\}.
\end{align*}
Note that the following map: 
$$ (d_1, d_2, d_3,d_4, \ldots,d_{2p-1}, d_{2p}) \longrightarrow  (-d_1, d_2,-d_3, d_4, \ldots, -d_{2p-1}, d_{2p}),$$
 is a bijection map from  $G^{(T)}_{2p,m}$ to $G^{(H)}_{2p,m}$. Hence from (\ref{eqn:card_D02pm_T})
	\begin{equation} \label{eqn:card_D02pm_H}
	\# G^{(H)}_{2p,m}= 	\sum_{\ell=0}^{\lfloor \frac{2mp}{2m+1} \rfloor} (-1)^\ell \binom{2p}{\ell} \binom{2p+2pm-(2m+1)\ell-1}{2p-1}. 
\end{equation}

Now  we prove Theorem \ref{thm:unbdd_LSD}.
\begin{proof}[Proof of Theorem \ref{thm:unbdd_LSD}]
	To show that a probability distribution $F$ has unbounded support, it is sufficient to show that $(m_{2p})^{1/p} \to \infty$ as $p \to \infty$, where $\{m_p\}$ is the moment sequence of $F$. 
	First note that $\# G^{(T)}_{2p,m}, \# G^{(H)}_{2p,m}  \geq 1$, where $\# G^{(T)}_{2p,m}$ and $\# G^{(H)}_{2p,m}$ are as in (\ref{eqn:card_D02pm_T}) and (\ref{eqn:card_D02pm_H}), respectively. Therefore from (\ref{eqn:moment_beta_T}) and (\ref{eqn:moment_beta_H}), we have $\beta^{(T)}_{2p}  \geq \gamma_T(p)$ and $\beta^{(H)}_{2p} \geq \gamma_H(p)$, where $\gamma_T(p)$ and $ \gamma_H(p)$ are as in (\ref{eqn:lim_2pmoment_T1}) and  (\ref{eqn:lim_2pmoment_H1}), respectively. Now using Proposition A.1 and Proposition A.2 from \cite{bryc_lsd_06}, we have that $(\gamma_T(p))^{1/p}$ and $ (\gamma_H(p))^{1/p}$ tend to $\infty$ as $p \to \infty$. These show that $(\beta^{(T)}_{2p})^{1/p}, (\beta^{(H)}_{2p})^{1/p}$ tends to $\infty$ as $p \to \infty$. This completes the proof of Theorem \ref{thm:unbdd_LSD}.
\end{proof}
\begin{remark}
	 Suppose symmetric Toeplitz and Hankel matrices have the entries $\{\frac{Y^{(m)}_j}{\sqrt n};{j\geq 0}\}$, where $Y^{(m)}_j= \sum_{r=-m}^{m} c_r X_{j+r}$, where $c_r$ are positive constants and $\{X_k\}$ are i.i.d. with mean zero and variance one. Then Theorem \ref{thm:toeESD} and Theorem \ref{thm:hanESD} can be generalized for such matrices, and the idea of the proofs will be similar.
\end{remark}

\noindent \textbf{Acknowledgment:} The author would like to express his heartfelt gratitude to Prof. Koushik Saha for reading the draft and providing valuable suggestions. The research work is supported by the fund:
NBHM Post-doctoral Fellowship (order no. 0204/10/(25)/2023/R$\&$D-II/2803).


\providecommand{\bysame}{\leavevmode\hbox to3em{\hrulefill}\thinspace}
\providecommand{\MR}{\relax\ifhmode\unskip\space\fi MR }
\providecommand{\MRhref}[2]{%
	\href{http://www.ams.org/mathscinet-getitem?mr=#1}{#2}
}
\providecommand{\href}[2]{#2}

\end{document}